# Optimal Load Restoration in Active Distribution Networks Complying with Starting Transients of Induction Motors

H. Sekhavatmanesh, *Member, IEEE*, J. Rodrigues, C. L. Moreira, J.A.P. Lopes, *Fellow, IEEE*, R. Cherkaoui, *Senior, IEEE*

*Abstract—* Large horsepower induction motors play a critical role as industrial drives in production facilities. The operational safety of distribution networks during the starting transients of these motor loads is a critical concern for the operators. In this paper, an analytical and convex optimization model is derived representing the starting transients of the induction motor in a semi-static fashion. This model is used to find the optimal energization sequence of different loads (static and motor loads) following an outage in a distribution network. The optimization problem includes the optimal control of the converter-based DGs and autotransformers that are used for the induction motor starting. These models together with the semi-static model of the induction motor are integrated into a relaxed power flow formulation resulting in a Mixed-Integer Second Order Cone Programming (SOCP) problem. This formulation represents the transient operational limits that are imposed by different protection devices both in the motor side and network side. The functionality of the proposed optimization problem is evaluated in the case of a large-scale test study and under different simulation scenarios. The feasibility and accuracy of the optimization results are validated using I) off-line time-domain simulations, and II) a Power Hardware-In-the-Loop experiment.

## NOMENCLATURE

### A. Parameters (unless mentioned, all are in p.u.)

| | |
|---|---|
| $D_i$ | Importance factor of the load at bus $i$. |
| $F_{ij}^{th}$ | Square of the nominal thermal ampacity limit of line $ij$. |
| $H_m$ | Inertial constant of the motor $m$ (sec.) |
| $k_{max}$ | Number of steps assigned to each starting motor |
| $kp_i(kq_i)$ | Active (Reactive) load voltage sensitivity at bus $i$. |
| $P_{i,t}^0(Q_{i,t}^0)$ | Active (Reactive) nominal load at bus $i$ at time $t$. |
| $f_{max,i}^{DG}$ | Ampacity limit of the DG converter at node $i$. |
| $Kd_m$ | Coefficient of the friction and windage loss in the motor $m$. |
| $L_{i,t}^0$ | Load energization status at node $i$ and time t according to the steady-state analysis (1/0). |
| $r_{ij}(x_{ij})$ | Resistance (Reactance) of line $ij$. |
| $S_{k,m}$ | Slip value of the motor $m$ at step $k$. |
| $\Delta S_m$ | Slip interval between two successive steps for the motor $m$. |
| $w_{re}, w_{op}$ | Weighting factors of the objective function terms. |

### B. Variables (unless mentioned, all are in p.u.)

| | |
|---|---|
| $F_{ij,k,m}$ | Square of current flow magnitude in line $ij$, at step $k$ during the starting of the motor $m$. |
| $Fp_{i,m}^{DG}$ ($Fq_{i,m}^{DG}$) | Square of the active/reactive current references of the DG converter at node $i$ for the starting of the motor $m$. |
| $L_{i,t}$ | Binary decision variable indicating if the load at node $i$ is energized at time t or not (1/0) |
| $P_{i,k,m}^D$ ($Q_{i,k,m}^D$) | Active (Reactive) load power at bus $i$, at step $k$ during the starting of the motor $m$. |
| $p_{ij,k,m}$ ($q_{ij,k,m}$) | Active (Reactive) power flowing in line $ij$, at step $k$ during the starting of the motor $m$. |
| $P_{i,k,m}^{Sub}$ ($Q_{i,k,m}^{Sub}$) | Active (Reactive) power from the substation node $i$, at step $k$ during the starting of the motor $m$. |
| $P_{i,k,m}^{DG}$ ($Q_{i,k,m}^{DG}$) | Active (Reactive) power injection from the DG node $i$, at step $k$ during the starting of the motor $m$. |
| $\Delta r_m$ | Integer variable representing the tap position of the autotransformer at the motor $m$. |
| $\tilde{t}_{k,m}$ | An approximation for the acceleration time of the motor $m$ until step $k$ (sec.). |
| $\widetilde{\Delta t}_{k,m}$ | An approximation for the time lengths of the step $k$ in the acceleration period of motor $m$ (sec.). |
| $T_{k,m}^{ele}$ | Electrical torque of the motor $m$, at step $k$. |
| $T_{k,m}^{mec}$ | Mechanical torque of the motor $m$, at step $k$. |
| $U_{i,k,m}$ | Square of voltage magnitude at bus $i$, at step $k$ during the starting of the motor $m$. |

### C. Indices

| | |
|---|---|
| $i,j$ | Index of nodes |
| $ij$ | Index of branches |
| $k$ | Index of slip step |
| $m$ | Index of nodes hosting motor loads in the off-outage area |
| $t$ | Index of time step |

### D. Sets

| | |
|---|---|
| $N$ | Set of nodes |
| $N_m$ | Set of nodes hosting motor loads |
| $N_p$ | Set of protected nodes |
| $N_s$ | Set of nodes hosting static loads |
| $N^*$ | Set of nodes in the off-outage area |
| $N_m^*$ | Set of nodes hosting motor loads in the off-outage area |
| $T$ | set of all the time samples in the restorative period |
| $W$ | Set of lines |
| $W_p$ | Set of protected lines |

## I. INTRODUCTION

Large-scale induction motors have been used for more than a century as industrial drives for compressors, pumps, fans, or blowers[1]. The large and reactive current driven by induction motors during the starting period can impose high risks both in the motor side and in the network side [2].

The analysis of the motor acceleration transients has been studied in the literature in the context of many power system applications. For example, the effect of motor reacceleration is studied in power system restoration [3]–[5], voltage and frequency stability of islanded microgrids [6]–[10], and protection settings within industrial facilities [1], [11]. In these papers, the behavior of the induction motor is evaluated using only time-domain simulations. These approaches would be very time-consuming if they are used for decision-making problems with huge and complex solution spaces. For such problems, it is needed to formulate the motor starting dynamics in an analytical way such that it can be integrated to an optimization problem.

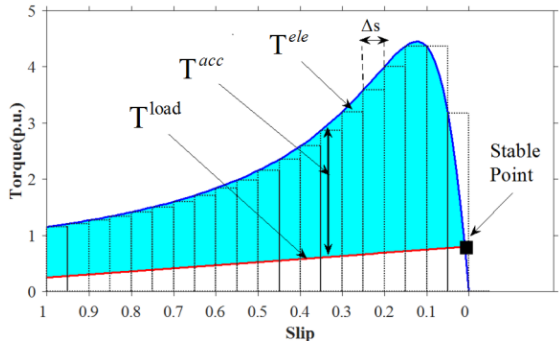

Fig. 1. The discretized toque-speed curve of the induction motor

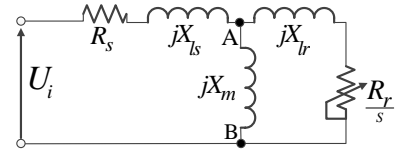

Fig. 2. The equivalent circuit of the induction motor.

Such analytical approaches are proposed in [12]–[15]. The authors of [12] developed nonlinear differential algebraic formulations to estimate the voltage dip during the motor starting. This voltage dip is predicted in [13] using neural network for an induction motor with a certain kVA capacity installed on a bus with a certain short circuit capacity. The developed formulations are all nonlinear and non-convex. Therefore, they cannot be integrated into the convex optimization problems. The presented formulation in [14] is incorporated into a maximum restorable load problem and solved iteratively using a heuristic approach. This approach is applicable only in the case of simple problems involving only one decision variable. The authors in [15] presented a quadratic optimization problem for the minimization of the voltage deviation with respect to the nominal value during the motor starting. In that paper, the motor reactive power during the acceleration period is approximated as a simple algebraic function of the terminal voltage magnitude. This simplified model does not account for the dynamics of the motor starting. In consequence, it cannot evaluate correctly the feasibility of the motor starting with respect to the operational safety constraints.

In this paper, we study the induction motor starting within the context of the load restoration problem. When a fault occurs in a distribution network, once it is isolated, the area downstream to the fault place remains unsupplied which is called the *off-outage area*. This area is re-energized by the healthy neighboring feeders using switching operations. This new configuration remains for a so-called restorative period until the faulted element is repaired. During this period, it is aimed to restore the loads in the most optimal way such that the total energy not supplied is minimized. In this regard, the network security constraints must be respected especially in case of starting large motor loads. The problem of finding the optimal load energization sequence given a new network configuration is referred as the load restoration problem.

Apart from the load energization sequence, another decision variable that can support the network security constraints during the motor acceleration is the control of converter-interfaced generation units. In [16]–[19], different control strategies are proposed for the control of the DG converter under grid fault condition. The aim is to support the low-voltage right through (LVRT) capability of the DG. For achieving this aim, these control strategies focus either only on the active power or only on the reactive power injection of the DG. This weakness is addressed in [20] and [21] by considering the resistive characteristics of the distribution line impedances and controlling both active and reactive power set points of the DG converter. However, these control methodologies are aimed to support the voltage only at a single node. Therefore, the network is simplified using its Thevenin's equivalent seen from that node.

In this paper, an analytical optimization model is derived for the load restoration problem. The main goal is to derive the optimal energization sequence of different type of loads. The problem is formulated as MISOCP and solved using Gurobi solver. Compared with the state-of-the-art algorithms proposed so far for the restoration problem, the major contributions of this paper are the following:

- The proposed load restoration problem incorporates a convex semi-static model of induction motor loads representing their starting transients.

- A convex model is derived for the converter-interfaced DGs working in constant current mode following voltage sags induced by motor load startings. This model is integrated into the optimal load restoration problem in order to obtain the optimal current set points of these DGs. Using these set points, the DGs support the electrical safety constraints in an optimal fashion during motor starting transients.

- The developed optimization problem includes a convex model of the autotransformer that is used for the starting of the induction motor. Therefore, the optimal tap setting of this autotransformer is derived from the proposed optimization problem.

- The transient operational limits imposed by under-voltage and over-current protection devices are integrated into the developed optimization problem. The aim is to guarantee that the starting transients of motor loads will not trigger the protection devices that exist in the distribution network.

## II.  MODELLING OF THE INDUCTION MOTOR STARTING

In this section, a semi-static model is proposed for the starting dynamics of the induction motor such that it can be integrated into the power flow formulation in a convex fashion. The aim is to formulate the operational safety constraints in the whole distribution network during the motor starting period. Fig. 1 shows the general trajectory of the electrical torque ($T^{ele}$) generated by an induction motor during its acceleration assuming that the motor voltage is fixed. The load torque ($T^{load}$) and the acceleration torques ($T^{acc}$) are also shown in this figure. The former is defined as the summation of the mechanical torque ($T^{mec}$) on the shaft and the friction and windage loss ($Kd_m(1-S_{k,m})$). The latter is the difference between the electrical and load torques. The electrical torque of the motor load at node $m$ is formulated in (1) using the equivalent circuit of the induction machine shown in Fig. 2, where $R_s$ and $X_{ls}$ represent the resistance and reactance of the stator; $R_r$ and $X_{lr}$ are the resistance and reactance of the rotor; and $X_m$ is the magnetization reactance [12].

$$T_{k,m}^{ele} = \frac{R_r U_{th}/S_{k,m}}{(R_r/S_{k,m} + R_{th})^2 + (X_{lr} + X_{th})^2} \quad (1)$$

where,

$S_{k,m} = (n_s - n_m)/n_s$

$U_{th} = \frac{U_i X_m^2}{R_s^2 + (X_{ls}^2 + X_m^2)}$

$R_{th} = \frac{R_s X_m^2}{R_s^2 + (X_{ls}^2 + X_m^2)}$

$X_{th} = \frac{R_s^2 X_m + X_{ls} X_m (X_{ls} + X_m)}{R_s^2 + (X_{ls}^2 + X_m^2)}$

$n_s$ and $n_m$ are the synchronous speed and the rotor speed of the motor at node $m$, respectively. $R_s$ and $X_{ls}$ represent the resistance and reactance of the stator. $R_r$ and $X_{lr}$ are the resistance and reactance of the rotor. $X_m$ is the magnetization reactance. $U_{th}$, $R_{th}$ and $X_{th}$ are the Thevenin's voltage square, resistance and reactance seen from the rotor terminals (AB in Fig. 2), respectively. Inserting the expressions of these equivalent parameters in (1) results in a formulation for $T_{k,m}^{ele}$ that is a linear function of the square of the voltage terminal ($U_i$) and non-linear function of the slip ($S$). As illustrated in Fig. 1, in order to make a convex model, we discretize the slip range from the standstill ($s = 1$) until the stable point into $K_{max}$ fixed steps with equal length of $\Delta s$. In this regard, the slip value during each step $k$ is assumed fixed and equal to a given value $S_{k,m}$, which is considered as a parameter. This assumption is justified according to the discussion provided in [2]. In this paper, the terminology of *slip step* (or shortly *step*) refers to each of these discretized intervals within the transient acceleration period of each energized motor load[1]. As discussed in [2], the step length $\Delta s$ should be small enough depending on the total inertia of the induction motor. The next step is to derive the time duration of a given step $k$ indicated by $\Delta t_{k,m}$. The inverse of this time length is obtained using the dynamic motion equation given in (2) [22].

$$\frac{1}{\Delta t_{k,m}} = \frac{-1}{2H_m \cdot \Delta S_m}\left(T_{k,m}^{ele} - T_{k,m}^{mec} - Kd_m(1-S_{k,m})\right) \quad (2)$$

$$\tilde{t}_{k,m} = \sum_{k^*=1}^{k} \widetilde{\Delta t}_{k^*,m} \quad (3)$$

Since the acceleration torque during each step is fixed, the time derivative of the slip is represented by $\frac{\Delta s}{\Delta t}$. In order to derive $\Delta t_k$, the piece-wise linear approximation method is used as explained in Appendix A. In this regard, we add (13)-(16) to the set of constraints, replacing $x$ and $f(x)$ with $\frac{1}{\Delta t_{k,m}}$ and $\Delta t_{k,m}$, respectively. Therefore, $\widetilde{\Delta t}_{k,m}$ is obtained in a linear way as an approximation for $\Delta t_{k,m}$.

In (3), we obtain an approximation for the acceleration time of motor $m$ until step $k$ ($\tilde{t}_{k,m}$) by adding the approximated time lengths of all previous steps to step $k$ ($\widetilde{\Delta t}_{k,m}$). The obtained variable $\tilde{t}_{k,m}$ will be used in section III.C to derive the transient voltage and current limits.

During each single step, since the slip and therefore all the parameters of the motor equivalent circuit shown in Fig. 1 are fixed, the electrical state variables can be represented in the phasor domain. The aim of the next section is to obtain the values of these state variables for a given step $k$ and in the whole distribution network using AC power flow equations. In this paper, we neglect the DC term imbedded in the starting current of the induction motor. This assumption is justified because of the low X/R ratio in distribution networks. Therefore, the DC term of the starting current disappears shortly.

III. RESTORATION PROBLEM FORMULATION

In this section, the load restoration problem is presented as an example to show how the proposed semi-static model of the motor starting could be integrated into the power flow formulation. Due to the difficulties of integrating the transient constraints related to the motor load starting in the full restoration problem formulation, a two-stage approach is exploited.

In the first stage, we solve the restoration problem according to the formulation presented in [23]. This restoration problem contains the model of passive and active elements only in steady-state conditions. Solving this optimization problem provides I) the optimal configuration of the network (line switching variables) and II) the optimal load restoration sequences during the restorative period. This stage of analysis is referred in this paper as the *steady-state analysis*.

In the second stage, which is studied in this paper, we take the restoration solution obtained from the first stage and we modify it concerning the transient constraints of motor load starting. In this stage, we modify only the obtained load restoration sequence, while considering the starting dynamics of induction motors, which were neglected in the first stage. It means that the line switching variables are fixed to the ones obtained from the steady-state analysis. We assume that considering the starting dynamics of motor loads does not affect the optimal network configuration that was obtained during the steady-state analysis. In this second stage, we just change the energization sequence of the loads in the most optimal way such that the transient constraints are respected during the starting of each induction motor in the off-outage area. The main objective is to minimize the resulting increase of energy not supplied referring to its value obtained from the steady-state analysis.

The main decision variables of the proposed optimization problem are three folds, namely, I) the energization sequence of different loads during the time ($L_{i,t}$), II) the optimal tap setting of the autotransformer that is used for the starting of the induction motor ($\Delta r_m$), and III) the optimal active/reactive current injections by dispatchable DGs ($Fp_{i,m}^{DG}/Fq_{i,m}^{DG}$). The optimization problem is formulated in the form of a MISOCP. The optimization problem in the first stage (steady-state analysis) includes power flow formulation for each time step $t$ during the restorative period. However, in the second stage (studied in this paper), the optimization problem includes the power flow formulation for each step $k$ in the acceleration period of each motor load $m$ in the off-outage area.

---

[1] These slip steps ($k$) should not be confused with time steps ($t$), which refer to low-resolution intervals within the whole restorative period.

In this paper, it is assumed that only one motor load can be started at time and only once the starting transient of any other motor disappeared. As suggested in [1], the motor loads in an industrial plant are categorized into groups mainly based on their functional processes. These groups of motor loads are considered to be restored in successive time steps with certain intervals. Only the motor loads that are in the same group are restored simultaneously. In this regard, a given motor load in the proposed optimization problem can represent a group of motor loads in the LV network. The dynamic parameters of this aggregated motor load are specified according to the strategy given in [8].

According the assumption mention above, to each motor load in the off-outage area a specific set of steps k ($k \in \{1,2,\ldots,k_{max}\}$) is assigned. Therefore, all the electrical state variables are indexed with $k$ and $m$. The constraints involving these indices should hold for all the steps and all the motors in the off-outage area.

Minimize:
$$F^{obj} = W_{re}.F^{re} + W_{op}.F^{op} \quad (4)$$

$$F^{re} = \sum_{i \in N}\sum_{t \in T} D_i.(L_{i,t}^0 - L_{i,t}).P_{i,t}^0 \quad (5)$$

$$F^{op} = \sum_{m \in N_m^*}\sum_{k=1}^{k_{max}}\sum_{ij \in W} r_{ij}.F_{ij,k,m} \quad (6)$$

Subject to:
$$\begin{cases} L_{i,t} \leq L_{i,t}^0, \; L_{i,t} \leq L_{i,t+1} : & i \in N^* \\ L_{i,t} = 1 & : i \in N\setminus N^* \end{cases} \forall t \in T \quad (7)$$

**Load Modeling** (8)
**AC power flow formulation** (9)
**DG Modeling** (10)
**Autotransformer Modeling** (11)
**Transient Constraints** (12)

The objective function ($F^{obj}$) is formulated in (4) as the weighted sum of the reliability ($F^{re}$) and operational ($F^{op}$) objective terms. As mentioned earlier, the main objective is to minimize the unsupplied energy of loads due to the shifting in their energization times, which is referred in this paper as the reliability objective term. This energy is calculated in (5) summing the power of the all the loads that are not restored ($L_{i,t} = 0$), whereas they were commanded to be restored according to the steady state analysis ($L_{i,t}^0 = 1$). Unlike the reliability term, the operational term has a very small weighting coefficient. This term is expressed in (6) as the total active line power losses in the distribution network. This term is included in the objective function just to satisfy the exactness condition stated in [24] for the relaxed AC power flow formulation. According to this condition, the objective function of a minimization (/maximization) problem should strictly increase (/decrease) with the total line power losses in the network. Therefore, the squared current variables will be bounded at optimal values, ensuring the exactness of the optimal solution.

Constraint (7) enforces the load at node $i$ and time $t$ to remain unrestored ($L_{i,t} = 0$) if it is commended to be so according to the results of the steady-state analysis ($L_{i,t}^0 = 0$). Moreover, according to (7), once a load is restored ($L_{i,t} = 1$), it should remain supplied during the rest of the restorative period ($L_{i,t+1} = 1$). The loads that are not in the off-outage area should remain always supplied.

### A. Load Modeling

Equations (8.a) and (8.b) extract, respectively, $P_{i,m}^0$ and $Q_{i,m}^0$, as the nominal active and reactive powers of the load at node $i$ and at the time instant $t$, when the motor $m$ is started ($L_{m,t} - L_{m,t-1} = 1$). At this time $t$, if the load at node $i$ is not restored ($L_{i,t} = 0$), then $P_{i,m}^0$ and $Q_{i,m}^0$ will be zero. The product of binary variables in (8.a) and (8.b) introduce non-linear terms. There terms are linearized according to the reformulation technique proposed in [25].

$$P_{i,m}^0 = \sum_{t \in T} P_{i,t}^0.L_{i,t}.(L_{m,t} - L_{m,t-1}) \quad \forall i \in N, \forall m \in N_m \quad (8.a)$$

$$Q_{i,m}^0 = \sum_{t \in T} Q_{i,t}^0.L_{i,t}.(L_{m,t} - L_{m,t-1}) \quad \forall i \in N, \forall m \in N_m \quad (8.b)$$

In the following, we use $P_{i,m}^0$ and $Q_{i,m}^0$ to formulate the active and reactive power consumptions of different type of loads. We start with the static loads. The active and reactive power of the static load at node $i$ and at step $k$, during the starting of motor load $m$ are expressed in (8.c) and (8.d), respectively. In this regard, the exponential model is used. Assuming that $U_{i,k,m}$ is close to $1\,p.u$, the model is linearized using the binomial approximation approach as proposed in [23]. The products of the binary variables $L_{i,t}$ (in the formulation of $P_{i,m}^0$ and $Q_{i,m}^0$) and the positive continuous variable $U_{i,k,m}$ introduce non-linear terms in (8.c) and (8.d). In order to preserve the linearity, these terms are re-formulated as given in Appendix B.

$$P_{i,k,m}^D = P_{i,m}^0\left(1 + (U_{i,k,m} - 1)\right)^{kp_i/2} \approx$$
$$\approx P_{i,m}^0\left(1 + \frac{kp_i}{2}(U_{i,k,m} - 1)\right) \quad \forall i \in N_s, \forall m \in N_m, \forall k \quad (8.c)$$

$$Q_{i,k,m}^D = Q_{i,m}^0\left(1 + (U_{i,k,m} - 1)\right)^{kq_i/2} \approx$$
$$\approx Q_{i,m}^0\left(1 + \frac{kq_i}{2}(U_{i,k,m} - 1)\right) \quad \forall i \in N_s, \forall m \in N_m, \forall k \quad (8.d)$$

Now, we move to the formulation of the motor load powers. The active and reactive power of the motor load at node $i$ and at step $k$, during the starting of motor load $m$ are expressed in the following. For the motor load that is starting ($i = m$), (8.e) and (8.g) express the active and reactive power consumptions according to the equivalent circuit shown in Fig. 2. The other motor loads ($i \neq m$) that are already energized are modeled in

(8.f) and (8.h) as PQ constant loads during the starting of the motor load $m$.

$$P_{i,k,m}^D = \begin{cases} U_{i,k,m}\left(\dfrac{R_{k,m}^{th}}{R_{k,m}^{th\ 2} + X_{k,m}^{th\ 2}}\right) & : i = m \\ P_{i,m}^0 & : i \neq m \end{cases} \quad \forall i, m \in N_m \quad \begin{matrix}(8.e)\\ \forall k \\ (8.f)\end{matrix}$$

$$Q_{i,k,m}^D = \begin{cases} U_{i,k,m}\left(\dfrac{X_{k,m}^{th}}{R_{k,m}^{th\ 2} + X_{k,m}^{th\ 2}}\right) & : i = m \\ Q_{i,m}^0 & : i \neq m \end{cases} \quad \forall i, m \in N_m \quad \begin{matrix}(8.g)\\ \forall k \\ (8.h)\end{matrix}$$

where, $R_{k,m}^{th}$ and $X_{k,m}^{th}$ represent the Thevenin's equivalent resistance and reactance seen from the terminals of motor $m$ at a given step $k$, respectively. The value of these Thevenin's equivalent impedances depend on the slip value at each step $k$.

### B. AC power flow formulation

Constraints (9.a)-(9.d) represent the second-order cone relaxation of the branch flow model proposed in [26] for each step $k$ in the acceleration period of each motor load $m$. The aim is to extract the electrical state variables (i.e. voltage, current, and power flow variables). These variables are needed for the optimal control of the converter-interfaced generation units and for checking the transient constraints as will be discussed in sections III.C and III.E, respectively.

$$U_{j,k,m} = U_{i,k,m} - 2(r_{ij} \cdot p_{ij,k,m} + x_{ij} \cdot q_{ij,k,m}) + F_{ij,k,m}(r_{ij}^2 + x_{ij}^2) \quad (9.a)$$
$$\forall ij \in W, \forall m \in N_m, \forall k$$

$$p_{ij,k,m} = \sum_{i^* \neq i} p_{ji^*,k,m} + r_{ij} \cdot F_{ij,k,m} + P_{j,k,m}^D - P_{j,k,m}^{Sub} - P_{j,k,m}^{DG} \quad (9.b)$$
$$\forall ij \in W, \forall m \in N_m, \forall k$$

$$q_{ij,k,m} = \sum_{i^* \neq i} q_{ji^*,k,m} + x_{ij} \cdot F_{ij,k,m} + Q_{j,k,m}^D - Q_{j,k,m}^{Sub} - Q_{j,k,m}^{DG} \quad (9.c)$$
$$\forall ij \in W, \forall m \in N_m, \forall k$$

$$0 \leq U_{i,k,m} \leq v_{max}^2 \quad \forall i \in N, \forall m \in N_m, \forall k \quad (9.e)$$

$$F_{ij,k,m} \geq \dfrac{p_{ij,k,m}^2 + Q_{ij,k,m}^2}{U_{i,k,m}} \quad \forall ij \in W, \forall m \in N_m, \forall k \quad (9.d)$$

Constraint (9.a) expresses the nodal voltage equation as given in [26]. The last term in the right hand side of (9.a) is usually neglected, since it is much smaller than the other terms. Constraints (9.b) and (9.c), respectively, concern with the active and reactive power balances at the extremities of each line. The first term in the right hand side of (9.b) and (9.c), represent, respectively, the sum of active and reactive power flows in lines that are connected to bus $j$ except the line $ij$. According to these two equations, the power flow from bus $i$ to bus $j$ is equal to the sum of power flows in other lines connected to bus $j$ plus the net injection power from bus $j$ and the power loss in line $ij$. Constraint (9.e) imposes the maximum voltage limit at the nodes of the network. Constraint (9.d) is the relaxed version of the current flow equation in each line according to [26]. This constraint is implemented in the form of the following second order cone constraint.

$$\left\| \begin{matrix} 2p_{ij,k,m} \\ 2q_{ij,k,m} \\ F_{ij,k,m} - U_{i,k,m} \end{matrix} \right\|_2 \leq F_{ij,k,m} + U_{i,k,m} \quad \begin{matrix}\forall ij \in W, \\ \forall m \in N_m, \forall k\end{matrix}$$

### C. Converter control of generation units

In this paper, we consider only dispatchable generation units including DGs, storage systems, and/or static synchronous compensators that are interfaced with the grid via full-bridge power converters. In this paper, the term of DG is used to refer to all these generation units. The re-energization of the DGs after the fault is accounted for in the steady-state analysis. In the optimization problem studied in this paper, it is assumed that the DGs are already connected to the grid and they work in *normal state*. It means that the DG converter controls the active and reactive power injections (or active power and voltage) following pre-determined set points. When a motor load starts, the voltage at the DG hosting node drops. This voltage drop is detected at the DG hosting node and launches the proposed voltage support control scheme, referred in this paper as *current saturation mode*. In this mode, instead of the active and reactive powers, the injection current of the DG is controlled. This control mode is achieved through already existing Fault-Right-Through (FRT) strategies in the converter interface while assuming that their current set points can be modified. According to the following formulation, we compute the optimal values of current set points in case of starting each of the motor loads during the restorative period.

$$Fp_{i,m}^{DG} + Fq_{i,m}^{DG} = f_{max,i}^{DG\ 2} \quad (10.a)$$

$$P_{i,k,m}^{DG\ 2} \leq Fp_{i,m}^{DG} \cdot U_{i,k,m} \quad (10.b)$$

$$Q_{i,k,m}^{DG\ 2} \leq Fq_{i,m}^{DG} \cdot U_{i,k,m} \quad (10.c)$$

The magnitude of the current injection by the DG converter at node $i$ is forced in (10.a) to be equal to its maximum current limit ($f_{max,i}^{DG}$). Constraints (10.b) and (10.c) derive the square of the active and reactive current components, respectively. These constraints are relaxed versions of the original formulations that are equalities instead of inequalities. The aim is to build a convex model of the DG control in current saturation mode. These relaxations will be exact according to the discussion provided in Appendix A. Constraints (10.b) and (10.c) are implemented in the form of second order cone constraints as expressed in the following:

$$\left\| \begin{matrix} 2P_{i,k,m}^{DG} \\ Fp_{i,m}^{DG} - V_{i,k,m} \end{matrix} \right\|_2 \leq Fp_{i,m}^{DG} + V_{i,k,m}$$

$$\left\| \begin{matrix} 2Q_{i,k,m}^{DG} \\ Fq_{i,m}^{DG} - V_{i,k,m} \end{matrix} \right\|_2 \leq Fq_{i,m}^{DG} + V_{i,k,m}$$

### D. Autotransformer tap setting

In some cases, the applied voltage to the induction motor terminal is reduced during the starting period using an autotransformer. This leads to reduce the starting current, power loss and radiated heat during the acceleration period. However, while reducing the voltage terminal, the starting torque will be reduced as well. It causes to lengthen the acceleration period. Therefore, there is a trade-off in setting the tap position of the autotransformer. In what follows, the autotransformer is modeled and incorporated into the optimization problem. The tap position of the autotransformer will be set according to the obtained optimal solution before starting the motor.

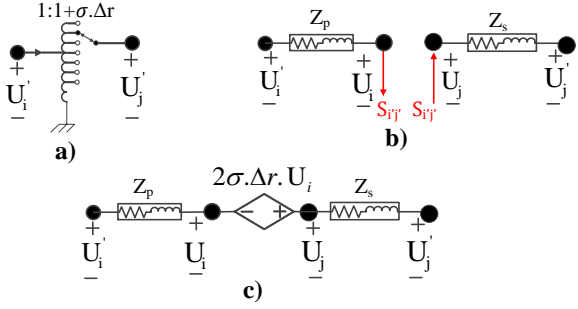

Fig. 3. Modelling of the autotransformer. a) schematic, b) standard equivalent circuit, c) linearized model.

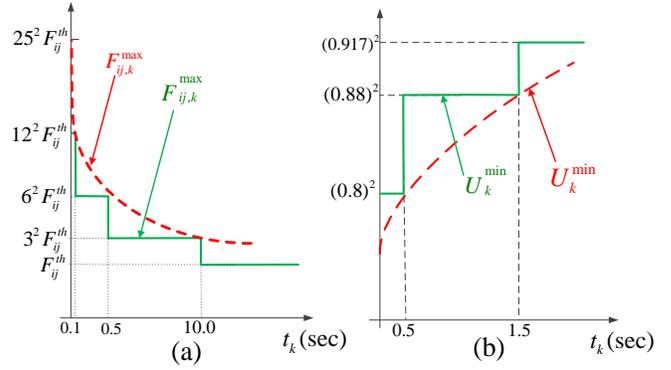

Fig. 4. Typical protection curve of a) over-current and b) under-voltage relays used in this paper

Fig. 3 shows the equivalent circuit of an autotransformer, where $Z_p$ and $Z_s$ denote the impedances at the primary and secondary sides, respectively. $S_{i'j'}$ is the apparent power flowing from the auxiliary node $i'$ to the auxiliary node $j'$. Based on this equivalent circuit, the constraints related to each autotransformer installed at the terminals of motor $m$ are formulated as follows:

$$U_{j,k,m} = U_{i,k,m} \cdot (1 + \sigma_m \cdot \Delta r_m)^2 \quad (11.a)$$

$$U_{j,k,m} \approx U_{i,k,m} \cdot (1 + 2\sigma_m \cdot \Delta r_m) \quad (11.b)$$

where, $\sigma$ is the ratio change between two consecutive taps of the autotransformer. Assuming that the turn ratio of the autotransformer is close to 1, the non-linear term in (11.a) is linearized in (11.b) using the binomial approximation. The product of the integer variable $\Delta r_m$ and the positive continuous variable $U_{i,k,m}$ makes a non-linear term. This non-linear term is reformulated according to the linearization strategy given in [23]. For this aim, first the integer variable $\Delta r_m$ should be expressed as the weighted sum of auxiliary binary variables. Then, the product of each of these auxiliary binary variables and the positive continuous variable $U_{i,k,m}$ is linearized using the approach given in Appendix B.

Finally, Fig. 3.c shows how the derived model in (11.b) together with the primary and secondary impedances of the autotransformer are incorporated into the AC-power flow equations.

*E. Transient Constraints*

In this section, transient constraints are expressed regarding the safe starting of the induction motor in a distribution network.

$$T_{k,m}^{ele} \geq T_{k,m}^{mec} + Kd_m(1 - S_{k,m}) \quad (12.a)$$

$$\widetilde{U}_{k,m}^{min} \leq U_{i,k,m} \quad \forall i \in N_p \quad (12.b)$$

$$F_{ij,k,m} \leq \widetilde{F}_{ij,k,m}^{max} \quad \forall ij \in W_p \quad (12.c)$$

In order to avoid the induction motor to stall, (12.a) enforces the electrical torque to be larger than or equal to the load torque. For a given slip $S_{k,m}$, the electrical torque of the motor $m$ is obtained using (1) and the mechanical load torque is determined according to the torque-speed curve of the mechanical load. This curve is assumed to be given for a specific load on the shaft.

In order to avoid the tripping of the under voltage relays and the over-current relays, (12.b) and (12.c) are added for every step k of the starting of each motor load m. Fig. 4 shows the typical protection curves reported in [5] and [10]. These curves represent the values of the under-voltage ($U_{k,m}^{min}$) and over-current ($F_{ij,k,m}^{max}$) limits as functions of the acceleration time ($\tilde{t}_{k,m}$), which was formulated in (3). In order to preserve the linearity in terms of variable $\tilde{t}_{k,m}$, $U_k^{min}$ and $F_{ij,k}^{max}$ are approximated by $\widetilde{U}_k^{min}$ and $\widetilde{F}_{ij,k}^{max}$, respectively, according to the piecewise linear approximation method explained in Appendix A. In this regard, for obtaining $\widetilde{U}_k^{min}$, we add (13)-(16) to the set of constraints, replacing $x$ and $f(x)$ with $\tilde{t}_{k,m}$ and $U_k^{min}$, respectively. For deriving $\widetilde{F}_{ij,k}^{max}$, we augment the set of constraints by (13)-(16), replacing $x$ and $f(x)$ with $\tilde{t}_{k,m}$ and $F_{ij,k}^{max}$, respectively.

IV. NUMERICAL RESULTS

In this section, the functionality of the proposed optimization model is evaluated using a test distribution network shown in Fig. 5. This network is based on a 11.4kV distribution network in Taiwan. The base power and energy values are set to 1MW and 1 MWh, respectively. The details regarding the nodal and line data are given in [27].

Except the motor loads that are indicated in Fig. 5, the rest of loads are assumed as static loads. The nameplate power ratings of the induction motors at nodes {13}, {20, 29}, and {41} are equal to 805, 435, and 80.5 horsepower, respectively. The parameters of the equivalent circuit (see Fig. 2) of these motor loads are adopted from the real data given in [1]. The static loads are assumed to be of constant-impedance type. Therefore, their load-voltage sensitivity coefficients are set to $K_p = 2$ and $K_q = 2$. It is assumed that each node in the network shown in Fig. 5 is equipped with a load breaker. The load variation data along time is according to the practical data reported in [27]. In Fig. 5, the critical loads are identified with '*'. The priority factor ($D_i$) of these loads is equal to 10 and for the other loads is equal to 1.

It is assumed that a fault occurs on the substation 86-84 shown in Fig. 5. The faulted substation is isolated by opening its nearest breakers at both sides. The restorative period is assumed from 17:00 P.M. to 03:00 A.M [28]. As mentioned in section III, in order to obtain the optimal restoration strategy, first, the steady-state analysis is performed according to [23]. The reconfiguration results are shown in Fig. 5. The areas in the off-outage area that are colored with red, green, blue, and yellow are energized through closing tie-switches T1, T2, T4, and T7. The optimal load restoration sequence obtained from the steady-state analysis is depicted on a time axes in Fig. 6.a.

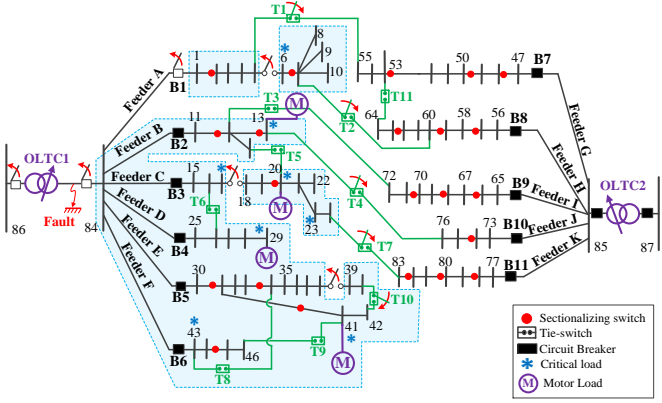

Fig. 5. The test distribution network under the post-fault configuration.

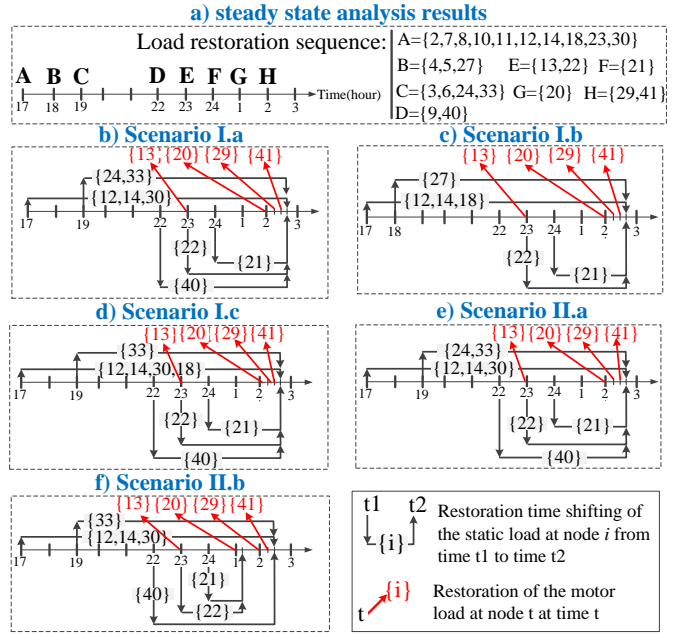

Fig. 6. The optimal load energization sequences.

Then, the proposed optimization problem in this paper is solved in order to shift the energization instants of certain loads to account for the starting dynamics of the motor loads. In this regard, three sets of simulation scenarios are studied. The algorithm is implemented on a PC with an Intel(R) Xeon(R) CPU and 6 GB RAM; and solved in Matlab/Yalmip environment, using Gurobi solver. Branch-and-Bound method is used to handle the developed mixed-integer optimization problem. The slip step size ($\Delta s$) is assumed to 0.05 p.u.

*A. Simulation scenario I: linear load toque*

The first set of simulation scenarios are studied assuming that the mechanical loads of all the motor loads have linear torque-speed characteristics. The magnitude of the nominal mechanical torques (at the synchronous speed) are equal to 0.4, 0.3, and 0.05 p.u. for the loads at nodes {13}, {20, 29}, and {41}.

Scenario I.a is defined assuming the critical and sensitive loads are set according to the default conditions shown in Fig. 5. The optimization problem is solved and the resulting optimal decisions are shown in Fig. 6.b. This figure shows the loads whose energization times should be shifted with reference to the results obtained from the steady-state analysis (see Fig. 6.a). The time resolution of this study is chosen to be 1 hour. As it can be seen in Fig. 6, the motor loads are energized in sequential steps with sufficient time intervals such that the starting transients of different motors do not overlap [1].

This shifting of the load energization times causes 29.7 p.u. additional energy not supplied (while considering the priority factors $D_i$). This simulation scenario includes a large-scale off-outage area and 4 unsupplied motor loads resulting in 134 binary and 24282 continuous variables. However, the solution for such a large study case is found just in 12.28 second.

In order to see the effects of the critical loads on the optimal results, scenario I.b is studied. In this scenario, the loads at nodes {24, 30, 33, 40} are considered as additional critical loads. The optimal results obtained for the load restoration sequence is depicted in Fig. 6.c. Compared to the results of scenario I.a (see Fig. 6.b), it can be seen that the energization times of the new critical loads are not postponed in scenario I.b. The reliability objective value is obtained equal to 29.95 p.u. The computation time is 11.64 second.

In the next step, scenario 1.c is defined such that node 18 is added to the protected nodes. The rest of the simulation conditions are the same as in scenario I.a. The results in Fig. 6.d shows that the restoration of the load at node 18 is shifted to the time when all the motor loads are already energized. The reason is that the under voltage limit at node 18 (during the motor starting transients) cannot be respected without shifting the energization of many of other loads. This decision leads to increase the reliability objective value from 29.07 p.u. in scenario I.a to 29.87 p.u. The computation time in scenario I.c is 11.69 second.

*B. Scenario II: fixed load torque*

In scenario II.a, it is assumed that the mechanical load on the shaft of the induction motor at node 20 has a fixed torque-speed characteristic equal to 0.06 p.u. The rest of the simulation conditions are the same as in scenario I.a. According to the results shown in Fig. 6.e, for the safe starting of the motor load at node 20, the same amount of loads should be shifted with respect to the ones in scenario I.a, although the mechanical load power in this scenario is so much less than the one in scenario I.a. The reason is that under a fixed-torque mechanical load, the induction motor can accelerate only if the starting torque (electrical torque at the standstill, s=1) is larger than the mechanical torque. In order to generate this starting toque, the voltage at the motor terminal should be large enough. This is obtained by shifting the energization of loads in the off-outage area, which results in a reliability objective value equal to 29.07 p.u.. The computation time is 16.98 second.

In Scenario II.b, the test case of the scenario II.a is studied while assuming an autotransformer at the terminals of the motor load at node 20. It is located in series between node 20 and the induction motor terminals during its acceleration period. Once the motor reaches to 80% of its nominal speed, this autotransformer is taken out from the circuit and the motor is directly connected to the grid [29]. This autotransformer enables $\pm 20\%$ voltage regulation range in 5 steps ($\sigma = 10\%, n = 4$). According to the results shown in Fig. 6.f, with the optimal setting of the autotransformer, there is no need any more to shift the energization time of loads at nodes {20, 24}.

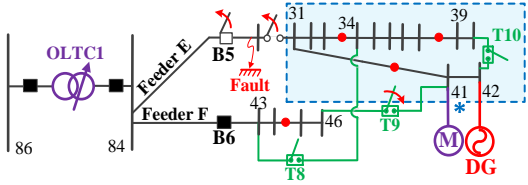

Fig. 7. Part of the test distribution network under the post-fault configuration in simulation scenario III.

Table I. The parameters of the induction motor at INESC TEC.

| Size (W) | $R_1(\Omega)$ | $X_1(\Omega)$ | $R_2(\Omega)$ | $X_2(\Omega)$ | $X_m(\Omega)$ | $H$(sec) |
|---|---|---|---|---|---|---|
| 4000 | 1.44 | 2.56 | 1.37 | 2.56 | 56.17 | 0.198 |

The optimal setting of the autotransformer is obtained with the tap position equal to -1. It means that the autotransformer reduces the voltage at motor terminals, which in turn reduces the starting current magnitude of the induction motor. Therefore, the magnitudes of the line voltage drops are reduced. In this regard, this optimal tap setting improves the quality of the restoration solution by increasing the margins of the nodal voltage magnitudes with respect to the transient under-voltage limits. The optimal value of the reliability objective is 25.175 p.u. and the computation time is 14.57 second. It will be validated in the next section that all the transient operational limits are respected with the obtained restoration solutions. Among all the constraints, the voltage magnitude at the accelerating motor node has a very narrow margin with respect to the transient minimum voltage limit.

*C. Scenario III: DG control in current saturation mode*

Scenario III is defined, where a fault occurs on line 84-30. Fig. 7 shows the part of the network that is affected by the post-fault configuration. As it can be seen, the only motor load in this part of the network is at node 41. The parameters of this induction motor are given in Table I. The mechanical load torque on the shaft of the induction motor changes linearly with speed and equals to 2.2 N.m. at the synchronous speed. There is also a DG at node 42, with the ampacity limit equal to 3.65 A. These parameters are according to the parameters of the physical induction motor and DG in the test setup at INESC TEC.

The proposed optimization problem is solved in case of this simulation scenario. The optimal load restoration sequence is to energize the static loads at nodes {32,35,36,37} together with the motor load at node {41} at the beginning of the restorative period (at 17:00 P.M). Once the transients of the motor starting disappear (at 17:15 P.M.), we will restore the loads at nodes {31,33,39,40}. The control of the DG converter at node 42 enters to the saturation mode during the motor acceleration period. The optimal active and reactive components of the converter current references (square of $Fp_{i,k}^{DG}$ and $Fq_{i,k}^{DG}$ used in (10)) are obtained as 2.106 and 2.970 ampere, respectively.

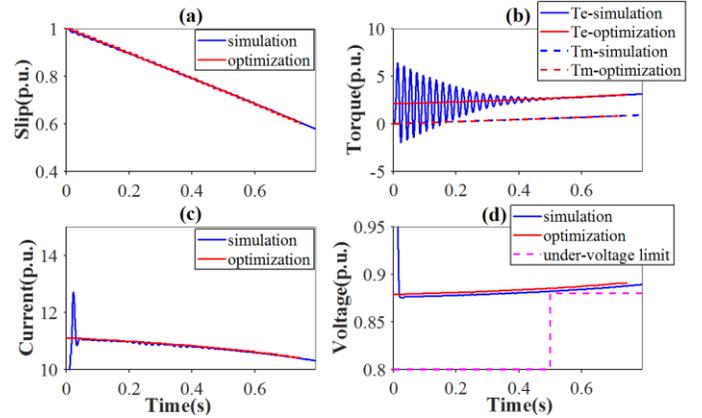

Fig. 8. Simulation results for scenario III. a) slip, b) electrical and mechanical torques, c) motor starting current, d) voltage at node 41.

## V. FEASIBILITY VALIDATION RESULTS

In this section, it is aimed to validate the solution feasibility of the proposed semi-static optimization model using an off-line simulation and a physical test experiment. It will be illustrated that the dynamics of the induction motor starting are represented into the optimization problem with a sufficient degree of accuracy. For this aim, the results obtained from scenario III and given in section IV.C are applied on the distribution network shown in Fig. 7.

*A. Time-domain simulation results*

In this section, an off-line model of the distribution network shown in Fig. 7 is built in Matlab/Simulink. The motor loads and static loads are represented by the Simulink model of the induction motor, and the impedance loads, respectively. The DG is modeled with a controllable dynamic load. First, we apply the obtained optimization results of scenario III on the off-line simulation model and then we start the motor. Fig. 8 shows the simulation results. In this figure, the electrical state profiles obtained from the optimization problem are compared with the ones obtained from the time-domain simulation. As it can be seen, the proposed semi-static optimization model does not represent the initial overshoot transients of the motor inrush current. As mentioned in section II, in deriving the semi-static model, we neglect the DC term in the motor starting current. Disregarding these very fast transients, Fig. 8 shows that the proposed semi-static model represents accurately the behavior of electrical state variables during the motor acceleration at the motor and network sides.

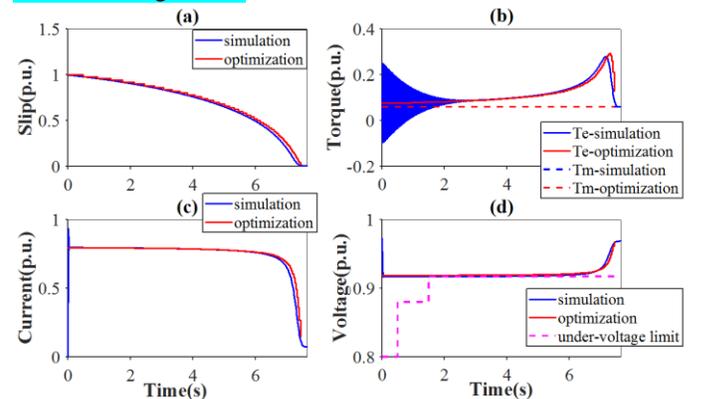

Fig. 9. Simulation results for scenario II.b. a) slip, b) electrical and mechanical torques, c) motor starting current, d) voltage at node 20.

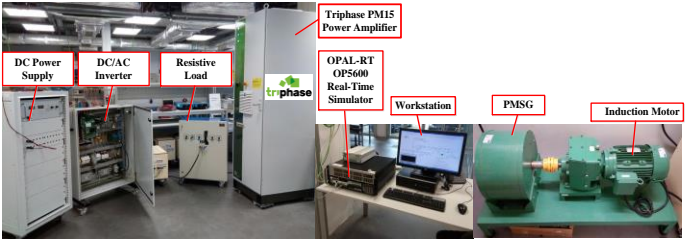

Fig. 10. The experimental PHIL test setup in the smart grid laboratory at INESC TEC.

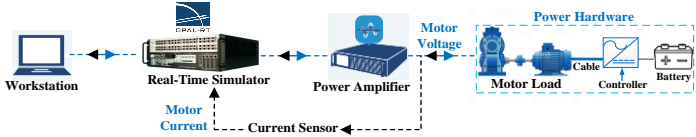

Fig. 11. Block diagram of the PHIL test setup.

In a similar fashion, the optimal results obtained from scenario II.b are tested using a time-domain simulation in Matlab/Simulink. The results at t=01:00, when the motor load at node 20 is started (see Fig. 6.f), are provided in Fig. 9. These results validate the feasibility of the optimal setting that was found in scenario II for the tap position of the auto-transformer.

*B. Experimental test*

As the next step of the validation study, a Power Hardware In the Loop (PHIL) experiment is performed. The test setup is implemented as shown in Fig. 10 in the smart grid laboratory at INESC TEC [30]. The block diagram of the laboratory test setup is depicted in Fig. 11. This PHIL test setup consists of three main parts. I) The first part includes the Matlab/Simulink network model of the whole network shown in Fig. 7 except the nodes 41 and 42. This model is compiled, and then executed by a Real-Time Simulator, namely OPAL-RT OP5600. II) The second part includes the Power Amplifier which magnifies the motor voltage signal received from the Real-Time Simulator. In this test setup TriPhase PM15 is used as the power amplifier with the nominal voltage equals to 400V. It can tolerate up to 30A of peak current per phase. In order to respect this ampacity limit during the motor starting transients, the base power and voltage of the test network are reduced to 320 VA and 150 V, respectively. III) The third part of the PHIL test setup is the physical hardware, including the motor load at node 41 and the DG at node 42 that are connected through a co-axial cable. In order to emulate a mechanical load with a linear torque-speed characteristic, the induction motor is loaded with a permanent magnet synchronous generator connected to a resistive load. The DG is emulated using an AC-DC inverter that is supplied by a DC power supply and controlled in current saturation mode.

We start the induction motor and measure the voltages at the motor terminals during the motor acceleration period. As shown in Fig. 12, this waveform is always above the transient under-voltage limit. It verifies that the under-voltage relay installed at the motor terminals will not trip in case of starting the motor load. In Fig. 13, the voltage magnitude measured at the motor terminal during the experiment is compared with the voltage profile obtained from the proposed semi-static optimization model. This figure shows that except at the first instants after the motor energization, the results of the semi-static model are sufficiently accurate with respect to the experimental results. At the first instants following the motor starting, the voltage at the motor terminal experiences a voltage dip that cannot be represented using the semi-static model. This voltage dip is caused by the DC term in the motor inrush current, which is neglected in deriving the semi-static model. As it can be seen in Fig. 13, this transient voltage dip disappears very fast and does not affect the feasibility of the solution with respect to the transient under-voltage limit.

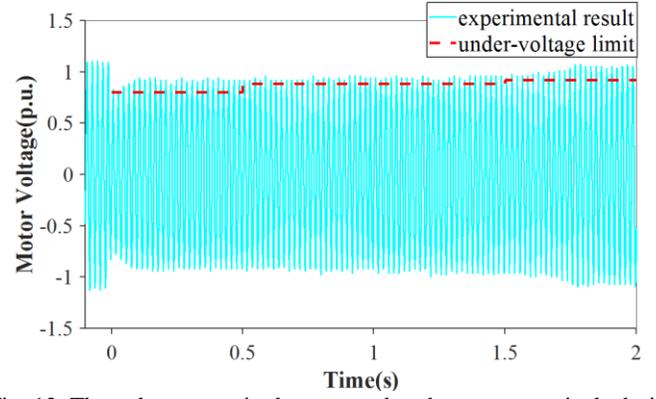

Fig. 12. The voltage magnitude measured at the motor terminals during its acceleration period.

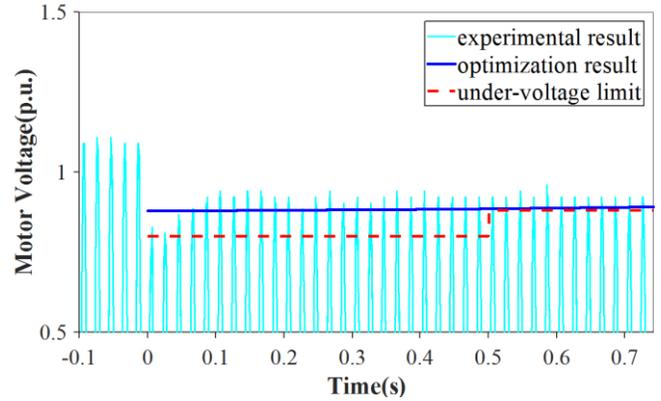

Fig. 13. The voltage at the motor terminal during its acceleration period obtained from the experiment and from the optimization model.

Fig. 14 shows the voltage, current, and power measured at the terminals of the DG during the motor acceleration period. The active and reactive components of the DG current are almost controlled to the reference values, which are obtained in section IV.C as 2.106 and 2.97 Ampere, respectively. Fig. 14 validates that the DG is working in current constant mode during the period that the DG experiences the voltage dip at its terminals due to the inrush starting current of the motor.

## VI. CONCLUSION

This paper presents a semi-static optimization model for the starting dynamics of induction motors. This model is integrated into the power flow formulation in a convex fashion. The resulting optimization problem represents the transient states of an active distribution network during the motor acceleration transients. The current saturation mode of the DG converter control is modeled in this formulation in a convex way. In this regard, the optimal references for the active and reactive components of the injection current are obtained. In addition, the proposed optimization formulation includes a convex model for the optimal tap setting of the autotransformer that is used for the starting of the induction motors.

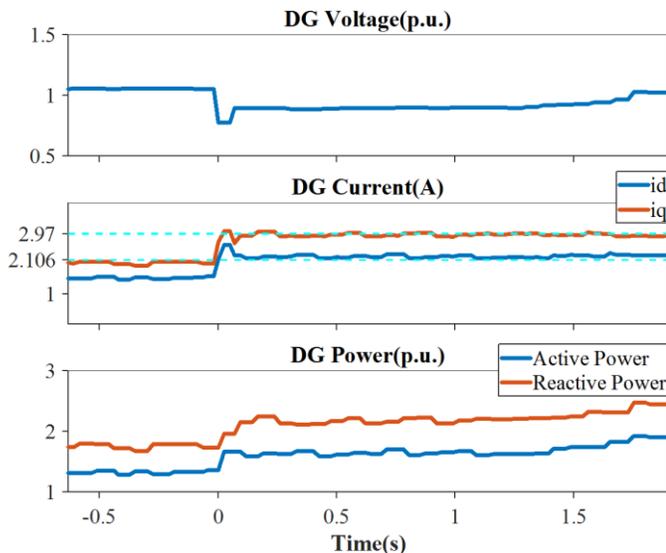

Fig. 14. The measured voltage, current and power of the DG during the motor acceleration period in the PHIL test experiment.

The resulting optimization formulation can be applied in any decision-making problem that is facing the dynamics imposed by the starting of induction motors. This model is used in this paper to solve the optimal load restoration problem in a distribution network. In this problem, the optimal energization sequence of different loads (static and motor loads) are obtained such that the total energy of loads that cannot be supplied will be minimized. The functionality of the proposed optimization model in solving large-scale load restoration problems is evaluated in the case of different simulation scenarios. The accuracy of the proposed semi-static optimization model in providing feasible solutions was verified using off-line time-domain simulations and also using a Power-Hardware-In-the-Loop test experiment.

## VIII. Appendices

### A. Piecewise Linear Approximation

Consider $f(x)$ as a continuous function with the domain of $[x_1, x_n]$. The concept of the piece-wise linear approximation introduced in [31] is shown in Fig. 15, where $\tilde{f}(x)$ denotes an approximation function for $f(x)$. We divide the domain of function $f$ by $n$ break-points $x_1, x_2, \ldots, x_n$. Between each two successive breaking points $x_i$ and $x_{i+1}$, the function $f$ is approximated with a straight line connecting the points $x_i$ and $x_{i+1}$. In order to formulate this approximation method, we introduce $n$ non-negative auxiliary variables $\lambda_i$ under (13)-(16).

Assume a given point $x$ between two breaking point $x_i$ and $x_{i+1}$. The value of $x$ is expressed in (13) in terms of the weighted sum of breaking points $x_i$ and $x_{i+1}$, with $\lambda_i$ and $\lambda_{i+1}$ as the weighting coefficients. Using the same auxiliary variables, function $\tilde{f}(x)$ is formulated as in (14). The sum of all auxiliary variables $\lambda_i$ should be one (15). Therefore, $\tilde{f}(x)$ consists of straight lines between successive breakpoints.

For the given point $x \in [x_i, x_{i+1}]$, all the auxiliary variables, except $\lambda_i$ and $\lambda_{i+1}$, will be zero. Because as given in (16), variables $\lambda$ are forced to respect the Special Ordered Set-2 (SOS2) constraint [32]. According to this constraint, out of all $\lambda_i$ variables, at most two successive variables can be non-zero. Most of the commercial solvers have the feature to account for these types of constraints during the Branch-and-Bound search algorithm.

$$x = \sum_{i=1}^{n} \lambda_i x_i \qquad (13)$$

$$\tilde{f}(x) = \sum_{i=1}^{n} \lambda_i f(x_i) \qquad (14)$$

$$\sum_{i=1}^{n} \lambda_i = 1 \qquad (15)$$

$$\lambda_i \geq 0 : SOS2 \qquad (16)$$

We can increase the accuracy of this approximation using more number of breaking points ($n$). However, it leads to more number of auxiliary variables, which could increase the computation burden of the optimization problem.

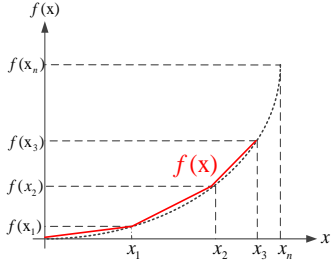

Fig. 15. Piece-wise linear approximation of an arbitrary continuous function $f(x)$

### B. Elimination of product of variables

In this section, a method is provided for the linearization of the constraints which incorporates a product of two variables. The product of two variables $x_1$ and $x_2$ can be replaced by one new variable y, subject to the constraints given in (17) and (18). The proof of this linearization is provided in [31]. It is assumed that $x_1$ is a binary variable and $x_2$ is a positive continuous variable, for which $0 \leq x_2 \leq u$ holds.

$$0 \leq y \leq ux_1 \qquad (17)$$
$$x_2 - u(1 - x_1) \leq y \leq x_2 \qquad (18)$$

### A. Exactness of the relaxation used in the modeling of the DG control in the current saturation mode.

According to [24], it can be proved that the relaxation made in (10.b) and (10.c) will be exact under the following conditions:

I) The maximum voltage limit at the DG hosting node is not binding during the acceleration period of the motor load.
II) The line ampacity limits are not binding during the acceleration period of the induction motor.
III) The objective function strictly decreases with the square of active and reactive power injections of the DG ($P_{i,k,m}^{DG}{}^2$ and $Q_{i,k,m}^{DG}{}^2$).

Condition I holds, since the voltage at the DG node is already dropped due to the starting of the induction motor. Therefore, there is a large margin for the voltage at the DG node with respect to its maximum limit. Condition II is usually ensured during the planning phase of the distribution networks. For example, in industrial networks, the line ampacity limits are sufficiently larger than the starting current of induction motors.

Under these two conditions, the optimal value of the objective function ($F_{re}$ in (4)) will be bounded only by the minimum voltage constraint. In this regard, during the starting of an induction motor, the minimum voltage limit at the hosting node of the starting motor is the bottleneck constraint. Therefore, in order to restore more loads (to decrease the value of objective function $F_{re}$), the voltage at the motor node should increase. From the other hand, with the starting of the motor load (increasing of power absorption at the motor node), the voltage at the DG terminal decreases. It means that the sensitivity of the voltage at the DG node ($V_{DG}$) is positive with respect to the power injection at the motor node ($P_m$: negative of power absorption). Due to the symmetry of the grid admittance matrix, we can intuitively infer that the sensitivity of the voltage at the motor node ($V_m$) with respect to the power injection at the DG node ($P_{DG}$) is also positive. The mathematical proof of this expression is beyond the scope of this paper.

Therefore, if the DG power injection increases, the voltage at the motor load increases and therefore the total unrestored energy of loads ($F_{re}$: objective function) decreases. It means that condition III also holds according to the specific characteristics of the load restoration problem studied in this paper.